\documentclass[12pt]{article}
\usepackage{latexsym}
\usepackage{amsfonts}
\def\dirac{{\slash \mkern-12mu D}} 

\def\ZZ{\mathbb Z}
\def\RR{\mathbb R}
\def\CC{\mathbb C}
\def\SS{\mathbb S}
\def\VV{\mathbb V}
\def\ro{\mathring{r}}
\def\ve{\varepsilon}
\def\bcp{\mathbb C \mathbb P}
\def\eea{\end{eqnarray*}}

\newtheorem{main}{Theorem}

\newtheorem{thm}{Theorem}[section]
\newtheorem{prop}[thm]{Proposition}
\newtheorem{cor}[thm]{Corollary}

\newenvironment{proof}{\medskip \noindent
{\bf Proof.}}{\hfill \rule{.5em}{1em}
\\}

\begin{document}
\sloppy
\title{Einstein Metrics,  Symplectic Minimality, \\
and Pseudo-Holomorphic Curves}

\author{Claude LeBrun\thanks{Supported 
in part by  NSF grant DMS-.} 
\\ 
SUNY Stony
 Brook 
  }

\date{}
\maketitle

\begin{abstract}
Let $(M^4, g, \omega )$ be a compact, almost-K\"ahler Einstein manifold
of negative 
star-scalar curvature.  Then $(M,\omega )$ is  a {\em minimal}   symplectic $4$-manifold
of general type. In particular, 
   $M$ cannot be differentiably decomposed as 
 a connected sum $N\# \overline{\bcp}_2$. 
\end{abstract}

\section{Introduction} 

A smooth Riemannian $n$-manifold $(M,g)$  is said to be {\em Einstein} if
its Ricci curvature, considered as a function on the unit tangent bundle of $M$,
is constant. If $n \geq 3$, this is equivalent \cite{bes} to saying that $g$ satisfies 
\begin{equation}
\label{albert}
r= \frac{s}{n}g,
\end{equation}
where $r$ and $s$ respectively denote the 
 Ricci tensor 
 and scalar curvature  of $g$, because equation 
  (\ref{albert}) and the contracted   Bianchi identity 
  $2 \mbox{ div } r = ds$
   together
 imply that   $s$ must be  constant if  $n\neq 2$.  As will be emphasized below, 
 the sign of this constant   plays an important  r\^ole  in much of 
  the theory of Einstein manifolds. 
 
 K\"ahler geometry provides  one of our richest sources of compact Einstein manifolds.
 Recall that a  Riemannian manifold is said to be {\em K\"ahler}
  if it admits a non-degenerate $2$-form
 with vanishing covariant derivative. When this happens, there must  in fact be 
 a non-degenerate closed $2$-form $\omega$ which can be expressed
 as $\omega=g(J\cdot , \cdot)$ for some integrable almost-complex structure
 $J$ on $M$.  One then says that the  closed $2$-form $\omega$ 
 is  the {\em K\"ahler form} of $(M,g,J)$, and its de Rham class
 $[\omega ] \in H^2 (M, \RR )$ is  called the {\em K\"ahler class}. 
 A celebrated result of Aubin \cite{aubin} and Yau \cite{yau} tells us precisely
 when a compact complex manifold $(M,J)$ admits a compatible K\"ahler-Einstein 
 metric  with $s < 0$; namely, this happens
if and only if  there is a holomorphic projective embedding $M\hookrightarrow \bcp_k$ 
for which the generic hyperplane section $M\cap \bcp_{k-1}$ is 
Poincar\'e dual to a negative multiple of $c_1(M, J)$.


When $n=4$, 
a K\"ahler-Einstein metric $g$ with $s < 0$ can  consequently exist   only  if   $(M,J)$ is 
a {\em minimal} complex surface of {\em general type}  \cite{bpv}. 
Recall that
 if $(N,J_N)$ is any complex
surface, and $p\in N$ is any point, 
 one may obtain a new complex surface $(\tilde{N}, J_{\tilde{N}})$
by replacing $p$ with a $\bcp_1$ of self-intersection $-1$. The resulting 
complex surface $\tilde{N}$ is called the {\em blow-up} of $N$ at $p$, and 
is diffeomorphic to the connected sum 
$N\# \overline{\bcp}_2$.
A complex surface is said to be {\em minimal} if it is {\em not} the blow-up of
some other complex surface.  Any compact complex surface $M$ can be obtained
from some minimal  complex surface $X$, called a {\em minimal model} for
$M$,   by  blowing up a finite number of times. 
A compact complex surface is said to be of {\em general type} if the dimension of 
the space of 
holomorphic sections of $K^{\otimes j}$ 
 is  quadratic in  $j \gg 0$,  where $K=\Lambda^{2,0}$
 denotes the {\em  canonical line bundle} of the complex surface. A key facet of   
Kodaira's classification theory \cite{bpv} is the 
assertion  that a   complex surface is  of general type iff it has a minimal 
model with   $c_1^2 > 0$ and with $c_1\cdot [\omega ] < 0$
 for some (and hence any) K\"ahler class.

 The purpose of the present article is to explore a
 symplectic generalization of the aforementioned link between 
complex-surface  minimality and the existence of K\"ahler-Einstein metrics. 
Recall that a  closed   non-degenerate $2$-form $\omega$ 
on a smooth manifold $M$ is called a 
 {\em symplectic form},
 and the pair  $(M,\omega )$ is 
then called a {\em symplectic manifold} \cite{mcsal}. 
The modern theory of $4$-dimensional symplectic manifolds
substantially parallels the theory of 
compact complex surfaces. In particular, 
if $(N,\omega_0 )$ is
a symplectic $4$-manifold, one can construct a 
$1$-parameter family of symplectic structures 
$\omega_{\epsilon}$ on 
the connect sum $\tilde{N}= N\# \overline{\bcp}_2$
 by removing a standard ball of volume $\epsilon^2/2$ 
 from  $(N,\omega_0 )$, and replacing it with a symplectic  
 $2$-sphere of area $\epsilon$ and self-intersection $-1$;
 any of the symplectic
manifolds $(\tilde{N}, \omega_{\epsilon })$ obtained by this construction 
is  then called \cite{mcrules,dusa} 
a {\em symplectic blow-up} of $(N,\omega )$. 
As in the complex case, every symplectic $4$-manifold can be 
obtained from a minimal model by blowing up a finite number of times, and,  
in keeping with the complex case,   a  symplectic
$4$-manifold is said  to be of {\em general type} \cite{lno,tjli} iff it has a minimal model
satisfying 
$c_1^2 > 0$ and $c_1\cdot [\omega ] < 0$. 

One of the most striking consequences of the existence of a symplectic
structure on a $4$-manifold is that it implies the existence of solutions of suitable perturbations
of the 
Seiberg-Witten equations, for any Riemannian metric \cite{taubes}. 
This in turn  leads to non-trivial constraints
on the Riemann curvature tensor \cite{lno,lcp2,lric,witten}, so the existence
of a symplectic structure exerts a ghostly influence over ostensibly unrelated
Riemannian geometries on a $4$-manifold. In particular, one obtains 
the following  obstruction to the existence of Einstein metrics \cite{lric}: 

\begin{thm} \label{corral}
Suppose that $(X, \omega_X)$ is a minimal symplectic $4$-manifold of general type.
Then the differentiable manifold $M = X \# k \overline{\bcp}_2$ does not
admit  Einstein metrics  if $k \geq c_1^2(X)/3$. 
\end{thm}

Thus,  the existence of Einstein metrics is obstructed 
on symplectic manifolds which are ``sufficiently non-minimal,'' 
even without requiring that the putative metric  be in any way
related to the symplectic form. In this article, we will see that one can  do 
distinctly better if 
 the metric and symplectic form are required to be related in 
a manner reminiscent of  K\"ahler geometry.

Any symplectic manifold $(M,\omega )$ admits almost-complex structures
$$J: TM\to TM, ~~ J^2=-{\bf 1},$$ such that the tensor field $g$ defined by 
\begin{equation}
\label{triple}
g (v, w) = \omega (v , Jw )
\end{equation}
is a Riemannian metric. One then says that the metric $g$ is {\em adapted} to 
$\omega$, and the triple $(M,g, \omega )$ is called 
an {\em almost-K\"ahler manifold}. If $(M,g)$ is a 
Riemannian $4$-manifold, then $(M,g, \omega )$ is almost-K\"ahler iff, 
for an appropriate orientation of $M$, 
$\omega$ is a  self-dual harmonic $2$-form with $|\omega | \equiv \sqrt{2}$.
If  $(M,g)$ is an oriented compact Riemannian $4$-manifold with $b_+(M) \neq 0$, 
the Hodge theorem therefore tells us that there is an open dense set
of $M$ on which $g$ is conformally related to an almost-K\"ahler metric,
since  the harmonicity and self-duality of $\omega$ are conformally invariant 
in dimension $4$. In other words, 
 the almost-K\"ahler condition 
  imposes no constraint at all on the local conformal geometry of  a  $4$-manifold. 
 It might  therefore seem implausible for such a weak condition to have a ponderable 
 Riemannian impact.

Nevertheless, 
 the so-called Goldberg  conjecture \cite{goldberg},
 that    
{\em every compact almost-K\"ahler Einstein manifold is  actually K\"ahler},
 remains completely open.  Moreover, 
Sekigawa \cite{seki} has shown that 
the  conjecture  is actually  true for metrics with 
 $s \geq 0$.  The $s< 0$ case remains a mystery, even in dimension $4$,
although some encouraging results have recently  been obtained concerning 
restricted  versions of the
problem in which the curvature tensor is also required to satisfy 
supplementary algebraic conditions  \cite{aad,arm1,arm2}.
On the other hand, it should be emphasized that, although  
Goldberg did not clearly stipulate that the manifolds in question 
are required to be compact, 
counter-examples  \cite{arm2}   due to K.P.  Tod 
show that the na\"ive local version of the conjecture 
is certainly false  in all dimensions $\geq 4$. 

If  one believed the $4$-dimensional case of the Goldberg conjecture to be true, 
one would consequently expect  almost-K\"ahler Einstein metrics 
with $s < 0$ to only exist on minimal complex surfaces. 
On the other hand,  if one believed the conjecture to be false, it might be 
helpful to  
know whether it would be a waste of time to look for counter-examples
on non-minimal symplectic $4$-manifolds. Thus, even for 
those who, like the author, remain agnostic as to the 
ultimate validity of the conjecture, it must seem reasonable to ask  
what current technology can tell us regarding  the  symplectic 
minimality of $4$-dimensional almost-K\"ahler
Einstein manifolds.  

The main  results of this paper involve the so-called {\em star-scalar
curvature} $s^\star$ of an almost-K\"ahler $4$-manifold
$(M,g,\omega )$,  as defined by
$$s^\star = 2 R (\omega , \omega ),$$
where $R$ is the Riemann curvature tensor of $g$. For a K\"ahler manifold,
this coincides with the usual scalar curvature 
$s$ of $g$, but for a  general almost-K\"ahler manifold one merely has $s^\star \geq s$. 
Our first main result is the following: 

\begin{main}\label{able}
Suppose that $(M,g,\omega)$ is a compact almost-K\"ahler Einstein $4$-manifold
with $s^\star < 0$. Then $(M,\omega)$ is a minimal symplectic $4$-manifold
of general type. Moreover, $M$ cannot be differentiably decomposed as
a connect sum $N\# \overline{\bcp}_2$. 
\end{main}

As is the case for Theorem \ref{corral},  the proof of  Theorem \ref{able} 
ultimately rests on a foundation of Seiberg-Witten theory, but in
the present case the  Seiberg-Witten equations enter 
only very  indirectly, via Taubes' existence theorem for
{\em pseudo-holomorphic curves} \cite{taubes3}. Let us now recall the 
definition of these last  objects, which 
were  first introduced into symplectic geometry 
by Gromov   \cite{gromsym}. Suppose that $(M,g, \omega )$ is an almost-K\"ahler manifold,
 and let $J$ be the associated almost-complex structure on $M$, 
 as determined by equation
(\ref{triple}). If $\Sigma$ is a compact  (but possibly disconnected) 
Riemann surface, with complex structure tensor $J_\Sigma : T\Sigma\to T\Sigma$,
one says that a  smooth, locally non-constant map
$\psi: \Sigma \to M$ is a pseudo-holomorphic curve
(or,  to emphasize the fixed choice of almost-complex structure, a  {\em $J$-holomorphic curve})
if its  differential $\psi_*$ is almost-complex-linear, in the sense that 
$$ J\circ  \psi_*  = \psi_*\circ J_\Sigma ~.$$
Taubes' results allow one to  predict the existence of a certain 
pseudo-holomorphic curve in any symplectic manifold  $M$
that can be decomposed as a connect sum $N\# \overline{\bcp}_2$, and Theorem \ref{able}
is then proved by calculating the evaluation of $c_1(M,J)$ on the 
homology class of this Riemann surface  in two different ways, with contradictory results.

In the absence of some {\em a priori} reason why 
a  compact 
almost-K\"ahler Einstein $4$-manifold with $s < 0$ should necessarily also have
$s^\star < 0$ everywhere, however, 
Theorem \ref{able} can be of only limited interest. 
Nonetheless, the same techniques used to prove 
Theorem \ref{able} 
lead to similar conclusions even when    the maximum of $s^\star$ is 
allowed to be  positive:

\begin{main}\label{baker} 
Suppose that  $(M,g,\omega)$ is an almost-K\"ahler  Einstein $4$-manifold 
with $s<0$ and  with 
$$ s^\star  \leq - \frac{s}{m}$$
for some constant $m >  1$.  Suppose, moreover, that $b_+(M)> 1$. 
Let $(X, \omega^\prime)$ be a 
symplectic minimal model for $M$, and let $k$ denote the number of 
blow-ups used to obtain $M$ from $X$.
Then 
$$k < \frac{1}{m+2}~c_1^2 (X)~.$$
\end{main}

Now  Oguri and Sekigawa \cite{ogsek}  claim
  that 
any almost-K\"ahler Einstein $4$-manifold with $s < 0$ automatically
satisfies $s^\star \leq -s/6$.
If this is correct, the hypotheses of Theorem \ref{baker}
should always  hold  for  $m=6$, thereby allowing one to conclude  
that blow-ups $X\# k \overline{\bcp}_2$ 
of symplectic manifolds $X^4$ with $b_+(X) > 1$ never   
  admit almost-K\"ahler Einstein metrics when $k\geq c_1^2(X)/8$. 
 This would certainly represent a considerable sharpening of Theorem \ref{corral} 
 in the almost-K\"ahler case. Better estimates for the global maximum of  
  $s^\star$ might  thus appear to represent a compelling topic  for future research.

\section{Almost-K\"ahler Geometry}
 
 Suppose that $(M,g,\omega )$ is an almost-K\"ahler 4-manifold.
 That is, we suppose that $g$ is a Riemannian metric, that $\omega$
 is  a symplectic form, and  that
  ${J_{a}}^{b}=\omega_{ac}g^{bc}$
 is an almost-complex structure. Orienting our $4$-manifold
 $M$ so that the $4$-form $\omega\wedge \omega$ is positive, 
 this is equivalent to saying that  $\omega$ is a closed self-dual 2-form
 on $(M,g)$ whose point-wise norm is everywhere $\sqrt{2}$.

The star-scalar curvature of $(M,g,\omega)$ is defined to be 
$$s^\star  = 2R(\omega , \omega ) = \frac{1}{2}R^{abcd}\omega_{ab}\omega_{cd},$$
where $R$ denotes the Riemann curvature tensor of $g$. 
This will coincide with the usual scalar curvature $s$
of $g$ if and only if  $(M,g,J)$ is a K\"ahler manifold. 
Indeed,  the Weitzenb\"ock formula for  a closed self-dual 2-form
tells us that  
$$0=\frac{1}{2}\Delta |\omega|^{2} + | \nabla \omega |^{2}
-2W^+ (\omega , \omega ) +\frac{s}{3} |\omega|^{2},$$
and plugging in $|\omega|^2\equiv 2$ thus yields 
$$s^\star  = 2\left[W^+(\omega , \omega ) + \frac{s}{12} |\omega|^2\right]= s + |\nabla \omega|^2,$$
so that $s^\star \geq s$, with equality everywhere precisely when 
$(M,g,\omega )$ is K\"ahler. 

The almost-complex structure $J$ induces a decomposition 
$$\Lambda^2_\CC = \Lambda^{1,1} \oplus \Lambda^{2,0} \oplus \Lambda^{0,2}$$
of the complex-valued $2$-forms, 
and  this is  related to the self-dual/anti-self-dual
decomposition
$$\Lambda^2 = \Lambda^+ \oplus \Lambda^-$$
induced by $g$ and the orientation by 
\begin{eqnarray*}
\Lambda^+_\CC & = & \CC \omega \oplus \Lambda^{2,0} \oplus \Lambda^{0,2}, \\
\Lambda^{1,1} & = &  \CC \omega \oplus \Lambda^-_\CC.
\end{eqnarray*}
In all these formul{\ae}, the summands are mutually orthogonal with respect to the 
Hermitian inner product induced by $g$. 
The line bundle $\Lambda^{2,0}$ is of particular interest, and is called 
the {\em canonical} line bundle $K$ of $(M,J)$.
Its dual $L=K^{-1}$ is called the {\em anti-canonical} line bundle, and 
may be naturally identified with $\overline{K}=\Lambda^{0,2}$
by using the  Hermitian inner product induced by $g$.

The anti-canonical line bundle $L$ of an almost-K\"ahler 4-manifold turns out to carry a natural 
connection, originally discovered by Blair \cite{blair}, and  later
rediscovered by Taubes \cite{taubes} via a  completely different line of reasoning.  
The purpose of this section is to point out a remarkable curvature property of 
this connection that holds when the metric $g$ is also Einstein. 

\begin{prop}\label{watson} 
Let $(M,g,\omega)$ be an almost-K\"ahler Einstein $4$-manifold.
Then the anti-canonical line bundle $L=K^{-1}$ of
$M$ carries a natural Hermitian connection whose curvature $2$-form 
$F$ satisfies  
$$iF = \eta + \frac{s+s^\star }{8}\omega + \frac{s-s^\star }{8}\hat{\omega},$$
 where 
 $$\eta = W^+(\omega )^\perp\in \Lambda^{2,0}\oplus \Lambda^{0,2},$$
while 
 $$\hat{\omega}\in \Lambda^-$$
 is  defined only on the open set $U=\{ p\in M~|~(s^\star - s)(p)\neq 0\}$,  satisfies 
 $$|\hat{\omega}|=|{\omega}|\equiv \sqrt{2},$$
 and corresponds to an almost-complex structure $\hat{J}$ on $U\subset M$
 which is 
 compatible with the {\em anti-symplectic} orientation of $M$. 
\end{prop}

We now give a concise proof of this result 
in terms of  the Penrose spinor calculus \cite{pr1,pr2,witten}, 
a brief summary of  which  may be found in 
 Appendix \ref{conv} below. 
 On any spin open subset of $M$,  we can uniquely express our self-dual harmonic 2-form as
 $$\omega_{ab}=\omega_{AB}\varepsilon_{A'B'},$$
 where
 $$\omega_{AB}=\bar{\omega}_{AB}=\omega_{(AB)}, ~~~~\omega^{AB}\omega_{AB}=2.$$
 Factoring $\omega$, we  therefore locally have
 $$\omega_{AB}= 2i\phi_{(A}\bar{\phi}_{B)},$$
where 
 $$\phi_{A}\bar{\phi}^{A} =1 .$$
 Now  $\omega$ is 
 co-closed, so 
 $$\nabla^{AA'}\omega_{AB}=0.$$
 But we also have 
 $$\nabla^{AA'}\varepsilon_{AB}=0,$$ 
and  $$\varepsilon_{AB}=2\phi_{[A}\bar{\phi}_{B]},$$
 so it follows that 
 $$0=\frac{1}{2}\nabla^{AA'}(-i\omega_{AB}+\varepsilon_{AB})=
 \nabla^{AA'}\phi_{A}\bar{\phi}_{B}.$$
 Contracting with $\phi^{B}$, we thus have
  $$0=\phi^{B}\nabla^{AA'}\phi_{A}\bar{\phi}_{B}
 =-\nabla^{AA'}\phi_{A}+\phi_{A}{\phi}^{B}\nabla^{AA'}\bar{\phi}_{B}.$$
 That is,
 $$(\nabla^{AA'}+\frac{i}{2}\vartheta^{AA'})\phi_{A}=0,$$
 where 
 $$\vartheta_{AA'}=2i{\phi}^{B}\nabla_{AA'}\bar{\phi}_{B}=
 2i\bar{\phi}^{B}\nabla_{AA'}\phi_{B}.$$
 Notice that $\vartheta=\bar{\vartheta}$, since $\overline{\bar{\phi}}_{A}=-\phi_{A}$. 
 We interpret this as saying that $\phi$  solves the 
Dirac equation $\dirac^\theta \phi =0$ 
 for a specific  spin$^{c}$ structure and a specific $U(1)$ connection $\theta$ on 
its determinant  line bundle $L=\Lambda^2{\mathbb V}_+$.
We can then view the local ambiguity $\phi_{A}\leadsto e^{iu/2}\phi_{A}$
 in choosing $\phi$ as simply stemming from a change of local 
 trivialization in the Hermitian  line bundle $L$,
 so that $\phi$ becomes a global unit section of 
 $${\mathbb V}_{+}={\mathbb S}_{+}\otimes L^{1/2},$$
 where $L$ is concretely the anti-canonical line bundle of 
 ${J_{a}}^{b}=\omega_{ac}g^{bc}$, and is endowed with a 
 standard connection $\theta$, with connection $1$-form   given by 
 $i\vartheta$ in our system of local trivializations.
 
 Now  
 $\nabla |\omega|^2=0$, so  
 $$\omega^{BC}\nabla_{AA'}\omega_{BC}=0,$$ 
 and hence 
 $$\phi^{B}\bar{\phi}^{C}\nabla_{AA'}\omega_{BC}=0.$$
 But since $\omega$ is co-closed, we have 
 $$\nabla_{AA'}\omega_{BC}= \nabla_{A'(A}\omega_{BC)},$$
and,  since $\phi$ and $\bar{\phi}$ form  a basis for  ${\mathbb S}_+$, 
 it  follows that \begin{equation}
\label{curious}
\nabla_{AA'}\omega_{BC}= \alpha_{A'}\phi_{A}\phi_{B}\phi_{C}
                      +\bar{\alpha}_{A'}
                      \bar{\phi}_{A}\bar{\phi}_{B}\bar{\phi}_{C}
\end{equation}
 for some unique $\alpha_{A'}$.
 Notice that the substitution $\phi \leadsto e^{iu/2}\phi$
 results in the transformation  $\alpha \leadsto e^{-3iu/2}\alpha$,
 so that $\alpha$ can invariantly be interpreted as a 
 section of $\SS_{-}\otimes L^{-3/2}$. In particular, 
 the anti-self-dual 2-form 
\begin{equation}
\label{mirror}
\hat{\omega}_{ab}= 
 \frac{2i}{|\alpha|^{2}}\alpha_{(A'}\bar{\alpha}_{B')}\ve_{AB}
\end{equation}
 induces an almost-complex structure $\hat{J}_a^{~b}=
 \hat{\omega}_{ac}g^{bc}$ on the open set $U\subset M$ 
 where $\nabla \omega\neq 0$. 
Since equation (\ref{curious}) implies  that   $|\nabla \omega |^{2}=
2|\alpha |^{2}$, it now follows that
 \begin{equation}
\label{size}
|\alpha|^{2}=\frac{s^{\star }-s}{2}.
\end{equation}

 The twisted spinor field $\alpha$ also completely 
 encodes the covariant derivative of $\phi$. Indeed,   
 $$ \phi^{B}\nabla_{AA'}\phi_{B} =-\phi^{B}\phi^{C}
 \nabla_{AA'}\phi_{B}\bar{\phi}_C=
 \frac{i}{2}\phi^{B}\phi^{C}\nabla_{AA'}\omega_{BC}=
 \frac{i}{2}\bar{\alpha}_{A'}\bar{\phi}_{A}$$
 and hence 
 \begin{equation}
 	\nabla_{AA'}\phi_{B} = 
 -\frac{i}{2}\left[ \bar{\alpha}_{A'}\bar{\phi}_{A}\bar{\phi}_{B}
 +\vartheta_{AA'}\phi_{B}\right] .
 	\label{der}
 \end{equation}
 In other words, 
 $$\hat{\nabla}_{AA'}\phi_{B}=  -\frac{i}{2}
 \bar{\alpha}_{A'}\bar{\phi}_{A}\bar{\phi}_{B} , $$
 where $\hat{\nabla}$  is the spin$^{c}$ connection on
 $\VV_{+}=\SS_{+}\otimes L^{1/2}$ induced by the $U(1)$ connection
 $\theta = i\vartheta$ on $L$.  
 
 Let us now calculate the curvature of $L$; cf. \cite{arm1,td}.
 Since $L$ carries the connection form $i\vartheta$, its
 curvature is given by $F=id\vartheta$, so that 
 \begin{eqnarray*}
 	iF^{+} & = & \varepsilon_{A'B'}\nabla^{C'}_{(A}\vartheta_{B)C'}  \\
 	iF^{-} & = & \varepsilon_{AB}\nabla^{C}_{(A'}\vartheta_{B')C} 
 \end{eqnarray*}
 Hence
 \begin{eqnarray*}
 	iF^{+}_{AB} & = &   \nabla^{A'}_{(A}\vartheta_{B)A'}\\
 	 & = &  \nabla^{A'}_{(A}2i\phi^{C}\nabla_{B)A'}\bar{\phi}_{C}  \\
 	 & = & -2i\phi^{C}\Box_{AB}\bar{\phi}_{C}
 +2i\left[\nabla^{A'}_{(A}\phi^{C}\right]\nabla_{B)A'}\bar{\phi}_{C}  \\
 	 & = & -2i \phi^{C} {W^+_{ABC}}^{D}\bar{\phi}_{D} -2i \phi^{C}\frac{s}{24} 
 (\varepsilon_{AC}\bar{\phi}_{B}+\varepsilon_{BC}\bar{\phi}_{A}) \\
 	 &  & +2i (-\frac{i}{2})\left[ 
 	 \bar{\alpha}^{A'}\bar{\phi}_{(A|}\bar{\phi}^{C}
 +\vartheta_{(A|}^{A'}\phi^{C} \right](\frac{i}{2})\left[
 \bar{\bar{\alpha}}_{A'}\bar{\bar{\phi}}_{|B)}\bar{\bar{\phi}}_{C}
 +\vartheta_{|B)A'}\bar{\phi}_{C}
 \right]   \\
 	 & = & {W^+_{AB}}^{CD}\omega_{CD}+\frac{s}{12} \omega_{AB}  
 	  - \frac{1}{4}|\alpha|^{2}\omega_{AB} 
 	  \\&=& \frac{s+s^{\star }}{8}\omega_{AB}+ \eta_{AB} 
 \end{eqnarray*}
 where $\eta = [W^+(\omega )]^\perp$ is orthogonal to $\omega$. 
 
Similarly, we  compute the anti-self-dual part of
the curvature of $L$:
  \begin{eqnarray*}
                     	iF^{-}_{A'B'}	 & = &  \nabla^{A}_{(A'}\vartheta_{B')A}  \\
                     	 & = &  \nabla^{A}_{(A'}2i\phi^{C}\nabla_{B')A}\bar{\phi}_{C}   \\
                     	 & = & -2i\phi^{C}\Box_{A'B'}\bar{\phi}_{C}
 +2i\left[\nabla^{A}_{(A'}\phi^{C}\right]\nabla_{B')A}\bar{\phi}_{C}    \\
                     	 & = &  i \phi^{C} {\ro_{A'B'C}}^{D}\bar{\phi}_{D}  \\
                     	 &  &  +2i (-\frac{i}{2})\left[ 
 	 \bar{\alpha}_{(A'|}\bar{\phi}^{A}\bar{\phi}^{C}
 +\vartheta_{(A'|}^{A}\phi^{C} \right](\frac{i}{2})\left[
 -{\alpha}_{|B')}{\phi}_{A}{\phi}_{C}
 +\vartheta_{|B')A}\bar{\phi}_{C}
 \right]  \\
                     	 & = & -\frac{1}{2}{\ro_{A'B'}}^{CD}(2i
                     	 \phi_{(C}\bar{\phi}_{D)}) 
                     	 -\frac{i}{2}\bar{\alpha}_{(A'}\alpha_{B')} \\
                     	 & = & -\frac{1}{2}{\ro_{A'B'}}^{AB}\omega_{AB}-
                     	  \frac{i}{2}{\alpha}_{(A'}\bar{\alpha}_{B')} 
                     \end{eqnarray*}  
  When $g$ is Einstein, $\mathring{r}=0$, and we thus have 
  $$iF^+= \eta+ \frac{s+s^\star }{8}\omega , ~~ iF^-=  \frac{s-s^\star }{8}\hat{\omega},$$
  and Proposition \ref{watson} follows.

 Another   remarkable  consequence  of these formul{\ae}  
  is the following \cite{blair}:
      
      \begin{prop}[Blair]
      Let $(M,g,\omega)$ be any $4$-dimensional compact
      almost-K\"ahler manifold. Then 
      $$
\int_M \frac{s+s^\star }{2}d\mu   = 4\pi c_1\cdot [\omega ]  ~.
$$
      \end{prop}
      \begin{proof}
      Either side can be rewritten as $2\int_M iF^{+}\wedge \omega$. 
      \end{proof}
      
      In particular, since $s^\star \geq s$, this implies 
      
\begin{cor} \label{blair}
For any compact
      almost-K\"ahler $4$-manifold $(M,g,\omega)$,
       $$
\int_M s^\star ~d\mu   \geq  4\pi c_1\cdot [\omega ]  \geq \int_M s ~d\mu ~.
$$    
\end{cor}

 The inequality $s^\star  \geq s$ 
is a repetitive motif  underlying much of almost-K\"ahler geometry, 
like a {\em basso ostinato}.
In light of our formulation of Proposition \ref{watson}, it would seem very natural to ask
whether one can ever actually have strict inequality $s^\star  > s$ everywhere. In
fact, this never happens in 
the Einstein case \cite{arm1}:

\begin{prop}[Armstrong] 
If $(M,g,\omega)$ is a compact $4$-dimensional almost-K\"ahler Einstein manifold,
there must be at least one point $p\in M$ at which  $s^\star =s$.
\end{prop}
\begin{proof}
Suppose not. Then $\alpha\neq 0$ on all of $M$, and
(\ref{mirror}) defines an 
 almost-complex structure 
$\hat{J}$ on the entire manifold. Moreover, since $\alpha\in {\mathbb S}_-\otimes L^{-3/2}$,
the anti-canonical line bundle of $\hat{J}$ is $L^{-3}$. 
Thus the reverse-oriented manifold $\overline{M}$ has an 
 almost-complex structure $\hat{J}$ with anti-canonical line bundle $L^{-3}$, whereas 
 $M$ has an almost-complex
 structure $J$ with anti-canonical line bundle $L$.
 Hence 
 $$(2\chi + 3\tau )(M)= \int_{M}c_{1}^{2}(L) =c_1^2(M),$$
 while 
 \begin{eqnarray*}
 (2\chi - 3\tau )(M)	 & = & (2\chi + 3\tau )(\overline{M})  \\
 	 & = & \int_{\overline{M}}c_{1}^{2}(L^{-3})  \\
 	 & = & -\int_{M}[-3c_{1}(L)]^{2}  \\
 	 & = & -9 \int_{M}c_{1}^{2}(L) = -9 c_1^2(M), ~
 \end{eqnarray*}
 so these two homotopy invariants of $M$ must have opposite signs. 
 But, on the other hand,   the Hitchin-Thorpe inequality  
 \cite{bes,hit} tells us that these characteristic numbers 
 must both be non-negative; indeed, since $g$ satisfies  $\mathring{r}=0$, the
 Gauss-Bonnet and signature formul{\ae} tell us that   
 $$(2\chi\pm 3\tau)(M) = \frac{1}{4\pi^2} \int_M \left(\frac{s^2}{24}+ 2|W^\pm|^2\right)d\mu~.$$
 Hence $(M,g)$ must be flat! 
 In particular, $s\equiv 0$ and $s^\star = 2W^+(\omega , \omega ) + \frac{s}{3}=0$,
 so $s^\star \equiv s$, in contradiction to our assumption.  
 \end{proof}
 
 We will specifically make use of the following  consequence: 

\begin{cor}
\label{strongarm}
On any   compact  almost-K\"ahler Einstein $4$-manifold, 
 there is some point at which 
$|W^+|^2\geq s^2/24$.  
\end{cor}
\begin{proof}
Since $2W^+(\omega , \omega ) =s^\star - \frac{s}{3}$, 
and since $|\omega |^2 = 2$, there is always an eigenvalue $\lambda$ of 
$W^+ : \Lambda^+ \to \Lambda^+$ with $|\lambda|\geq |3s^\star - s|/12$. 
Since $W^+$ is trace-free, this implies that 
$$|W^+|^2 \geq \frac{3}{2}\lambda^2 
\geq \frac{(3s^\star - s)^2}{96},$$
and  at any  point where $s^\star = s$ we therefore must have
$|W^+|^2 \geq s^2/24$. 
\end{proof}

\section{Pseudo-Holomorphic Curves}

We will now  use Proposition \ref{watson} to estimate
 the pull-back of  the curvature
of $L$ to  a pseudo-holomorphic curve.

\begin{prop}\label{sherlock}
Let $(M,g,\omega)$ be a $4$-dimensional 
almost-K\"ahler Einstein manifold, and let $\psi : \Sigma \to M$
be a $J$-holomorphic curve in $M$. Then the pull-back $\psi^*F$ of the curvature
 of the anti-canonical line bundle $L\to M$ satisfies
 $$\frac{s}{4}\psi^*\omega ~\leq ~i\psi^*F ~\leq~ \psi^* \frac{s^\star }{4}\omega$$
 at each point of the oriented $2$-manifold $\Sigma$. 
\end{prop}
\begin{proof}
By Proposition \ref{watson}, the curvature of $L$ satisfies 
$$iF = \eta + \frac{s+s^\star }{8}\omega + \frac{s-s^\star }{8}\hat{\omega},$$
where $\eta\in \Lambda^{2,0}\oplus \Lambda^{0,2}$,
and where the last term  is understood to mean zero at the locus  $s^\star = s$
 (even though the bounded  $2$-form $\hat{\omega}$ is undefined there). 
Since $\psi$ is assumed to be $J$-holomorphic, we therefore have 
 $\psi^*\eta \equiv 0$, and hence 
$$i\psi^* F =\psi^*\left[ \frac{s}{4}\left(\frac{\omega + \hat{\omega}}{2}\right)
+ \frac{s^\star }{4}\left(\frac{\omega - \hat{\omega}}{2}\right)\right],
$$
Now   $\omega \in \Lambda^+$ and  $\hat{\omega} \in \Lambda^-$
correspond to almost-complex structures $J$ and $\hat{J}$ via
index-lowering with $g$.  
Since   $\Sigma$ is $J$-holomorphic, $\psi^*\omega$ is exactly the 
area form of the (possibly degenerate) metric $\psi^*g$, and we therefore have  
$$    - \psi^*\omega \leq \psi^* \hat{\omega} \leq \psi^*\omega $$
by applying  Wirtinger's inequality 
\cite{GH,HL} 
to $\hat{\omega}$. 
We  may therefore write  
$$   \psi^* \hat{\omega} = (1-2t) \psi^*\omega$$
for some $t \in  [0,1] $ at any point where $s^\star \neq s$, 
and it therefore follows that  
$$ i\psi^*F =\frac{(1-t)s  + ts^\star }{4}    ~\psi^* \omega ,~~ \exists t\in [0,1], $$
at each point of $\Sigma$. 
The inequality $s\leq s^\star $ therefore  implies that 
$$\frac{s}{4} \psi^* \omega \leq i\psi^* F \leq \psi^*\frac{s^\star }{4}  \omega $$
everywhere, exactly as claimed. 
\end{proof}

\begin{cor}
Let $(M,g,\omega)$ be a compact $4$-dimensional almost-K\"ahler Einstein manifold, and
suppose that   $\psi : \Sigma \to M$ is a $J$-holomorphic curve. Let  
$[\Sigma ]$ denote, interchangeably,  either  the image of the fundamental cycle of $\Sigma$
in $H_2(M,\ZZ)$ or its Poincar\'e dual in $ H^2(M, \ZZ)$. 
Then 
$$s ~[\omega ]\cdot [\Sigma ] ~\leq ~8\pi ~ c_1 \cdot [\Sigma ] ~\leq ~\int_\Sigma  s^\star  \omega.$$
\end{cor}
\begin{proof}
By the Chern-Weil theorem, the closed $2$-form  $iF$ represents 
$2\pi c_1=2\pi c_1(L)$ in de Rham cohomology. Now integrate on $\Sigma$, and apply Proposition 
\ref{sherlock}. 
\end{proof}

In particular, this implies the following:

\begin{cor}\label{evidently}
Let $(M,g,\omega)$ be an almost-K\"ahler Einstein $4$-manifold
with $s^\star  \leq 0$, and let   $\psi : \Sigma \to M$ be any $J$-holomorphic curve.
Then 
$$c_1\cdot [\Sigma ] \leq 0.$$
\end{cor}

We will now invoke Taubes' existence results for pseudo-holomorphic curves  \cite{taubes3}. 
One of Taubes' most striking and 
fundamental technical  results may be stated as follows: 

\begin{thm}[Taubes]
Let $(M,g,\omega)$ be an almost-K\"ahler $4$-manifold,
and let $F_0$ denote the curvature of its anti-canonical 
line bundle $L$. Let $\ell \to M$ be some  Hermitian complex line bundle,
and consider the    spin$^c$-structure with determinant
line bundle $L\otimes \ell^2$ obtained by setting 
$\tilde{\mathbb V}_+ = {\mathbb V}_+\otimes \ell =  {\mathbb S}_+\otimes L^{1/2}\otimes \ell $. 
If, for this spin$^c$ structure, 
there exists a sequence 
$(\Phi_j, \theta_j)$ 
of solutions of the perturbed Seiberg-Witten equations
\begin{eqnarray*}
\dirac^{\theta_j}\Phi & = & 0 \\
iF^+_{\theta_j}+\sigma (\Phi) & = & iF^+_0 
+t_j\omega 
\end{eqnarray*}
for some  sequence of  real numbers    $t_j\to \infty$,
then there is a $J$-holomorphic curve $\psi: \Sigma \to M$ in $(M,g,\omega )$
such that $[\Sigma ] = c_1(\ell )$. 
\end{thm}
In fact,  Taubes shows that a subsequence of  $\{ i(F_{\theta_j}-F_0)/4\pi\}$  converges as
a current to the desired  pseudo-holomorphic curve. 

The consequences of this result are particularly clean and dramatic when $b_+(M)> 1$.
In this setting, each spin$^c$ structure has a well-defined number attached
to it, called its Seiberg-Witten invariant \cite{taubes3}, which  roughly speaking 
counts the number of solutions, modulo gauge equivalence, of the perturbed 
 Seiberg-Witten
equations 
\begin{eqnarray*}
\dirac^{\tilde{\theta}}\Phi & = & 0 \\
iF^+_{\tilde{\theta}}+\sigma (\Phi) & = & \varphi \end{eqnarray*}
associated with a generic self-dual $2$-form $\varphi$. When this
invariant is non-zero, the Chern class of the given spin$^c$ structure is called
a {\em Seiberg-Witten basic class} \cite{witten}, and the 
perturbed Seiberg-Witten equations  are certainly guaranteed to have solutions for the 
choices of $\varphi$ used by Taubes. Thus:

\begin{cor} \label{basic}
Let $(M,g,\omega)$ be a compact almost-K\"ahler $4$-manifold
with $b_+\geq 2$, and suppose that $c_1+2a\in H^2 (M, \ZZ)/\mbox{\rm torsion}$ is 
a Seiberg-Witten basic class of the oriented $4$-manifold $M$. 
Then there is a  $J$-holomorphic curve $\psi: \Sigma \to M$ in $(M,g,\omega )$
such that $[\Sigma ]= a$. 
\end{cor}

When $b_+(M)=1$, the situation is a bit  more subtle.  If one considers
only those perturbations $\varphi$ for which the   harmonic part $\varphi_H$ has 
huge norm, the same recipe used to define the Seiberg-Witten invariant 
in the previous case   works perfectly well in most respects; however,  the value 
of the answer now depends not only on the oriented smooth $4$-manifold $M$, but also on
the  nappe (connected component) ${\mathcal C}^+$ of the cone 
$${\mathcal C} = \{ h \in H^2 (M, \RR )~|~ h\cdot h > 0\}$$
which  contains the de Rham class $[\varphi_H]$. To apply
Taubes' work when $b_+(M)=1$, one thus consider the 
so-called perturbed Seiberg-Witten invariants \cite{lcp2,liliu1,liliu2} of   $(M, {\mathcal C}^+)$,
 where  ${\mathcal C}^+$ is the nappe of $\mathcal C$  
containing  the symplectic class $[\omega]\in H^2(M, \RR)$. 

Now let $(M,\omega)$ be any compact symplectic $4$-manifold, and 
consider the standard spin$^c$ structure on $M$ induced by $\omega$. 
Another  result of Taubes \cite{taubes} tells us that the Seiberg-Witten invariant
of $M$ (respectively, the perturbed Seiberg-Witten invariant  of   $(M, {\mathcal C}^+)$)
 is non-zero for  this spin$^c$  structure, provided that $b_+(M)\geq 2$ (respectively, 
 provided that  
$b_+(M)=1$). 
Now if $M$ can be expressed as a connected sum $M\approx N \# \overline{\bcp}_2$,
there are self-diffeomorphism $\Psi: M\to M$ which, away from the neck connecting the
two  summands,  is given by  the identity on $N$, and by  complex conjugation
 on $\overline{\bcp}_2$.
If $b_+(M)=b_+(N)=1$, such a diffeomorphism moreover sends ${\mathcal C}^+$ 
to itself, since the positive sector of $H^2(M)$ actually arises from $H^2(N)$.
Moving the standard spin$^c$ structure via $\Psi$ thus gives rise to 
a new spin$^c$ structure for which the relevant version of the Seiberg-Witten
invariant is non-zero, and  we thus  obtain the following \cite{liliu1,liliu2,taubes3}: 

\begin{cor}\label{except} 
 Let $(M,g,\omega)$ be a compact almost-K\"ahler $4$-manifold. If $M$ is diffeomorphic to a 
connected sum $N \# \overline{\bcp}_2$, then there is a
$J$-holomorphic curve $\psi: \Sigma \to M$ in $(M,g,\omega )$
 such that $c_1\cdot [\Sigma ] =+1$. 
\end{cor}

For a generic choice of $\omega$-compatible  metric $g$, 
 the  $J$-holomorphic curve $\Sigma$ 
of Corollary \ref{except} will 
 be an embedded symplectic $2$-sphere of 
self-intersection $-1$;  this then allows one to symplectically blow
down $M$ to obtain  a symplectic structure on $N$. 
A symplectic $4$-manifold $M$ is said to be minimal if it cannot be
obtained from another symplectic manifold $N$ by  blowing up \cite{dusa}.
 As we have just seen, $(M,\omega )$ is 
minimal iff it contains no symplectic $2$-sphere of self-intersection 
$-1$, and  this  occurs iff 
 $M$ cannot be differentiably 
 decomposed as a connected sum $N\# \overline{\bcp}_2$.
 Any compact symplectic $4$-manifold can be obtained from
 a minimal one by blowing up a finite number of times. 

For our purposes, 
it is important to emphasize that Corollaries \ref{basic} and \ref{except}
produce pseudo-holomorphic curves  for
a fixed almost-K\"ahler metric $g$, and not just  for its generic perturbations; the
price that must be paid for this is simply that   the images  
  $\psi (\Sigma )$ of these curves may in principle be rather singular. 
  With this in mind, we now easily obtain the following result:

\begin{thm} \label{minimal} 
Let $(M^4,g,\omega)$ be an almost-K\"ahler Einstein manifold with 
$s^\star \leq 0$. Then $(M, \omega)$ is symplectically minimal. 
\end{thm}
\begin{proof}
If $M$ were non-minimal, Corollary \ref{except} would predict the 
existence of a $J$-holomorphic curve $\psi : \Sigma \to M$
with  $c_1\cdot [\Sigma ]=+1$. 
But since we have assumed that $(M,  g, \omega )$ is 
almost-K\"ahler and Einstein, with $s^\star \leq 0$,  Corollary
\ref{evidently} says we must have $c_1\cdot [\Sigma ] \leq 0$, so  this is a contradiction. 
\end{proof}

A  symplectic $4$-manifold is said to be
of {\em general type} \cite{lno,tjli}
 if it has a minimal model with $c_1^2 > 0$ and $c_1\cdot [\omega ] < 0$. 
This definition is chosen so that it coincides with Kodaira's
notion of general type when $M$ admits a complex structure. 
Note that the 
the condition that $c_1\cdot [\omega ] < 0$ is actually redundant \cite{taubes2} when
$b_+(M)\geq 2$. 

\setcounter{main}{0}
\begin{main}
If $(M,g,\omega)$ is a $4$-dimensional  almost-K\"ahler Einstein manifold with 
$s^\star  < 0$, then $(M,\omega )$ is a minimal symplectic manifold of 
general type, and   $M$ cannot  be differentiably decomposed
as a connect sum $N\# \overline{\mathbb C \mathbb P}_2$. 
\end{main}
\begin{proof}
Since $s^\star < 0$, Corollary \ref{blair} tells us that $c_1\cdot [\omega ] < 0$.
Because $g$ is assumed to be Einstein, and since
the assumption that  $s^\star < 0$  implies that $g$ has  $s < 0$,  
 we also know that $c_1^2 = (2\chi + 3\tau )(M) > 0$
by the Hitchin-Thorpe inequality. The claim thus follows  from 
Theorem \ref{minimal}. 
\end{proof}

We now consider the case when $s^\star$ is allowed to have a positive maximum.

\begin{main}
Suppose that  $(M,g,\omega)$ is an almost-K\"ahler  Einstein $4$-manifold 
with $s<0$ and  with 
$$ s^\star  \leq - \frac{s}{m}$$
for some constant $m >  1$.  Suppose, moreover, that $b_+(M)> 1$. 
Let $(X, \omega^\prime)$ be a 
symplectic minimal model for $M$, and let $k$ denote the number of 
blow-ups used to obtain $M$ from $X$.
Then 
$$k < \frac{1}{m+2}~c_1^2 (X)~.$$
\end{main} 
\begin{proof}
By assumption, $M$ contains $k$ disjoint  symplectic $2$-spheres of self-intersection $-1$; we 
we will  use  $E_1, \ldots , E_k$ to denote,  interchangeably,   
either the homology classes $\in H_2 (M, \ZZ)$ of these $2$-spheres,  or their  Poincar\'e duals 
$\in H^2 (M, \ZZ)$. Let $D\in H^2 (M, \ZZ)$ denote the pull-back of $-c_1(X, \omega^\prime)$ 
via the blowing-down map $M\to X$. Then $\pm D \pm E_1 \pm \cdots \pm E_k$
are all Seiberg-Witten basic classes of $M$, for all possible choices of 
the various $\pm$ signs in this expression. Since the first Chern class of 
$(M,\omega )$ is given by $c_1 = -D - E_1 - \cdots - E_k$, it 
thus follows that 
 $c_1 + 2D$, $c_1+ 2E_1$, \ldots , $c_1+2E_k$ are all basic classes. 
 Corollary \ref{basic}  therefore tells us that $D$, $E_1$, \ldots , 
 $E_k$ are all represented by $J$-holomorphic curves in  $(M,g,\omega )$. 
 Since our assumptions include the point-wise estimate  $s\leq - ms^\star $, 
 Proposition \ref{sherlock} thus tells us that 
 \begin{eqnarray*}
c_1^2(X) & = &D\cdot D \\
 & = &  -c_1\cdot D \\
 &\leq& -\frac{s}{8\pi} [\omega ] \cdot D \\
 & = &  -\frac{s}{8\pi} [\omega ] \cdot  (  - c_1 - \sum_{j=1}^kE_j)\\
 &=& \frac{s}{8\pi}\int_M \frac{s+s^\star }{8\pi}d\mu +  \frac{1}{8\pi}\sum_{j=1}^k\int_{E_j} s\omega 
 \\ &\leq&  \frac{1}{64\pi^2}\int_M s^2d\mu + \frac{1}{64\pi^2}\int_M ss^\star d\mu -  \frac{m}{8\pi}\sum_{j=1}^k\int_{E_j} s^\star \omega 
  \\ &\leq&  \frac{1}{64\pi^2}\int_M s^2d\mu + \frac{1}{64\pi^2}\int_M ss^\star d\mu -  m
  \sum_{j=1}^kc_1\cdot E_j
   \\ &\leq&  \frac{1}{64\pi^2}\int_M s^2d\mu + \frac{1}{64\pi^2}\int_M ss^\star d\mu ~~~ -  mk ~.
\end{eqnarray*}
 However, Proposition \ref{watson} tells us that 
 \begin{eqnarray*}
 c_1^2(M) &=& \frac{1}{4\pi^2} \int_M \left(|F^+|^2 - |F^-|^2\right)d\mu 
 \\ &=& \frac{1}{4\pi^2} \int_M \left(|\eta|^2 + \left[\frac{s+s^\star }{8}\right]^2|\omega|^2- 
  \left[\frac{s-s^\star }{8}\right]^2|\hat{\omega}|^2\right)d\mu 
 \\ &=& \frac{1}{4\pi^2} \int_M \left(|\eta|^2 + \frac{ss^\star }{8}\right)d\mu 
 \\ &\geq & \frac{1}{32\pi^2} \int_M ss^\star  d\mu ~,
 \end{eqnarray*}
 while the Gauss-Bonnet and signature theorems tell us that 
  \begin{eqnarray*}
 c_1^2(M) &=& (2\chi + 3\tau ) (M) 
 \\ &=& \frac{1}{4\pi^2}\int_M \left(\frac{s^2}{24}+  2|W^+|^2 - \frac{|\mathring{r}|^2}{2} \right)d\mu 
  \\ &=& \frac{1}{4\pi^2}\int_M \left(\frac{s^2}{24}+  2|W^+|^2 \right)d\mu 
   \\ &> & \frac{1}{96\pi^2}\int_Ms^2 d\mu ~,
 \end{eqnarray*}
 where in the last step we have used Corollary  \ref{strongarm},
 which asserts that   our
 almost-K\"ahler Einstein manifold $s < 0$ cannot have  $W^+\equiv 0$.
 Hence 
  \begin{eqnarray*}
c_1^2(X) &\leq&  \frac{1}{64\pi^2}\int_M s^2d\mu + \frac{1}{64\pi^2}\int_M ss^\star d\mu -  mk 
\\ & <  & \frac{3}{2} c_1^2(M) + \frac{1}{2} c_1^2(M) -mk 
\\ &=& 2c_1^2 (M) -mk 
\\ &=& 2c_1^2 (X) -2k -mk ~, 
\end{eqnarray*}
and it follows that  
$$k < \frac{ c_1^2 (X)}{m+2} ~~,$$
exactly 
as claimed. 
\end{proof}

\appendix

\section{Appendix: Spinor Calculus}
\label{conv}

If $(M,g)$ is an oriented  Riemannian $4$-manifold, and if 
$U\subset M$ is contractible open subset,
then $U$ carries spinor bundles ${\mathbb S}_{\pm}$
such that 
$$\CC \otimes TU= {\mathbb S}_{+}\otimes {\mathbb S}_{-}.$$
These bundles have structure group $SU(2)$, and carry 
standard connections induced by the Levi-Civita connection of
$(M,g)$. Because the representation theory of $SU(2)$  is so 
remarkably simple, it is often convenient  to derive 
facts about Riemannian $4$-manifolds by  exploiting
these local spinors, even if  $M$ itself  is not spin. 
This is particularly true when, as 
in the present article, $M$ carries some natural  
  {\em spin$^{c}$} structure, so that one has   
 global bundles
${\mathbb V}_{\pm}\to M$ which are formally given by
$${\mathbb V}_{\pm}={\mathbb S}_{\pm}\otimes L^{1/2}$$
for some Hermitian line bundle $L$ with 
$$c_{1}(L)\equiv w_{2}(M)\bmod 2.$$

The spinor calculations in this paper  are  carried out using  a modified version
 of  Penrose's 
spinor calculus  \cite{pr1,pr2}.
 This notation makes  use of  
{\em abstract indices}, which are best 
understood as  simply ``addresses'' which  label  
 the various input ``mailboxes''
of a spinor or tensor field.  Upstairs lower-case Roman 
indices $a,b,c, \ldots$ 
refer to the tangent bundle,  
upper-case unprimed-primed capital 
 Roman 
indices $A,B, C, \ldots$ refer to the 
spin bundle ${\mathbb S}_{+}$, and 
upper-case unprimed-primed capital 
 Roman 
indices $A', B', C', \ldots$ refer to the 
spin bundle ${\mathbb S}_{-}$.
Downstairs indices refer to  the dual
of the relevant bundle. 
Square brackets ``[~~]'' indicate skew-symmetrization
$$\varphi_{[\underbrace{ab\cdots c}_n]}:=\frac{1}{n!}
\sum_{\sigma\in S_n}(-1)^{\sigma}\varphi_{\sigma(a)\sigma(b)
\cdots \sigma(c)}~ , $$ whereas 
round brackets ``(~~)''  indicate   symmetrization:
$$\varphi_{(\underbrace{ab\cdots c}_n)}:=\frac{1}{n!}
\sum_{\sigma\in S_n}\varphi_{\sigma(a)\sigma(b)
\cdots \sigma(c)}~.$$ 
The repetition of an index, once upstairs and once downstairs, 
indicates contraction. 
The fundamental isomorphism 
$\CC \otimes T= {\mathbb S}_{+}\otimes {\mathbb S}_{-}$
is invoked by the convention that lower case Roman indices
$a,b,c, \ldots$ are to be viewed as equal to  the corresponding pairs 
$AA', BB', CC', \ldots$  For convenience, we  allow ourselves
to shuffle primed indices through unprimed indices, 
so that $A'AB$, $AA'B$, and $ABA'$ are all
regarded as the same; however,  it should be emphasized that
these are  {\em a priori}  different from e.g. $A'BA$. 

Because the bundles ${\mathbb S}_{\pm}$
have 
 structure group  
$SU(2)= Sp(1)$, they may 
be viewed as quaternionic line bundles. However, we 
instead choose to view them as $2$-dimensional 
complex vector bundles  equipped with 
complex anti-linear maps
$$
\begin{array}{ccccccc}
	{\mathbb S}_{+} & \longrightarrow  & {\mathbb S}_{+} &  & 
	{\mathbb S}_{-} &  \longrightarrow & {\mathbb S}_{-}  \\
	\beta^{A} & \mapsto & \bar{\beta}^{A} & ~~~~~~~~~ & 
	\gamma^{A'} & \mapsto & \bar{\gamma}^{A'}
\end{array}
$$
 which satisfy 
\begin{equation}
\bar{\bar{\beta}}^{A}= - {\beta}^{A}, ~~ \bar{\bar{\gamma}}^{A'}=
-\gamma^{A'}~~.	
	\label{joe}
\end{equation}
and correspond to multiplication by the quaternion
$j$. Since  $SU(2)\subset SL(2, \CC )$, 
 the bundles $\Lambda^{2}{\mathbb S}_{\pm}^{*}$ also 
admit canonical parallel sections $\varepsilon_{AB}=\varepsilon_{[AB]}$
and  $\varepsilon_{A'B'}=\varepsilon_{[A'B']}$, which are 
real 
$$\varepsilon_{AB}= \bar{\varepsilon}_{AB}, ~~ 
\varepsilon_{A'B'}= \bar{\varepsilon}_{A'B'}$$
with respect to the induced real structure,
and which are related to the Riemannian metric $g$ by
$$
g_{ab}= \varepsilon_{AB}\varepsilon_{A'B'}~. 
$$
There are also inverse sections $\varepsilon^{AB}$
and  $\varepsilon^{A'B'}$
of $\Lambda^{2}{\mathbb S}_{\pm}$, defined so that 
$${\varepsilon_{A}}^{B}=\varepsilon_{AC}\varepsilon^{BC},
~~{\varepsilon_{A'}}^{B'}=\varepsilon_{A'C'}\varepsilon^{B'C'}$$
are the identity endomorphisms of ${\mathbb S}_{\pm}$. 
Spinor indices are raised and lowered according to the 
convention that
$$\beta_{B}= \beta^{A}\varepsilon_{AB}, ~~ \beta^{A}= 
\varepsilon^{AB}\beta_{B}, $$
etc. 
 The norm of 
the spinor $\beta^{A}$ is defined by
$$| \beta |^{2}= \beta_{A}\bar{\beta}^{A}= - 
\beta^{A}\bar{\beta}_{A} \geq 0,$$
and the norm of $\gamma^{A'}$ is defined analogously. 

Any two form $\varphi_{ab}=\varphi_{[ab]}$ 
can be uniquely expressed as 
$$\varphi_{ab}= \varphi^{+}_{AB}\varepsilon_{A'B'}+
\varphi^{-}_{A'B'}\varepsilon_{AB}, $$
where 
$\varphi^{+}_{AB}=\varphi^{+}_{(AB)}$
and 
$\varphi^{-}_{A'B'}=\varphi^{-}_{(A'B')}$.
With  our orientation conventions, 
the $2$-forms 
$$\varphi^{+}_{ab}= \varphi^{+}_{AB}\varepsilon_{A'B'}$$
and 
$$\varphi^{-}_{ab}= \varphi^{-}_{A'B'}\varepsilon_{AB}$$
are then respectively  the {\em self-dual} and {\em anti-self-dual}
parts of $\varphi$,  in accordance with the decomposition 
$$\Lambda^{2}= \Lambda^{+}\oplus \Lambda^{-}$$ 
of the $2$-forms into the   eigenspaces
of  Hodge  star operator $\star$. Also notice that 
the usual convention on the norm of a $2$-form simplifies
in a rather pleasant way: for a real self-dual $2$ form 
$\psi_{ab}$, one has
$$
|\psi |^{2}:=\frac{1}{2}\psi_{ab}{\psi}^{ab}=
\frac{1}{2}\psi_{AB}\varepsilon_{A'B'}{\psi}^{AB}\varepsilon^{A'B'}=
\psi_{AB}{\psi}^{AB},$$
since ${\varepsilon_{A'}}^{A'}=2$.
Of course, the expression for the norm of an anti-self-dual
$2$-form simplifies in a similar manner. 

In much the same way, the Riemann curvature tensor  
decomposes as 
\begin{eqnarray*}
 R_{abcd}	 & = & W^+_{ABCD}\varepsilon_{A'B'}\varepsilon_{C'D'}+
   W^-_{A'B'C'D'}\varepsilon_{AB}\varepsilon_{CD}\\
 	 &  & -\frac{1}{2}\ro_{ABC'D'}\varepsilon_{A'B'}\varepsilon_{CD}
 	   -\frac{1}{2}\ro_{A'B'CD}\varepsilon_{AB}\varepsilon_{C'D'} \\
 	 &  & +\frac{s}{12} (\varepsilon_{AC}\varepsilon_{BD}
 	 \varepsilon_{A'C'}\varepsilon_{B'D'} 
 	-\varepsilon_{AD}\varepsilon_{BC}\varepsilon_{A'D'}\varepsilon_{B'C'})
 \end{eqnarray*}
 where 
 $W^{+}_{ABCD}=W^{+}_{(ABCD)}$, 
 $W^{-}_{A'B'C'D'}=W^{-}_{(A'B'C'D')}$,
 and $\ro_{ABA'B'}=\ro_{(AB)(A'B')}$. 
 Here 
 $$W^{+}_{abcd}=W^+_{ABCD}\varepsilon_{A'B'}\varepsilon_{C'D'}$$
 is the self-dual Weyl curvature, 
 $$W^{-}_{abcd}=
W^-_{A'B'C'D'}\varepsilon_{AB}\varepsilon_{CD}$$
is the anti-self-dual Weyl curvature,
$$s={R^{ab}}_{ab}$$
is the scalar curvature, and 
$$\ro_{ab}=\ro_{ABA'B'}={R^{c}}_{acb}-\frac{s}{4}g_{ab}
$$
is the trace-free part of the Ricci curvature.
Note that $W^{+}_{abcd}$ is usually considered \cite{bes} 
as a section of $\Lambda^{+}\otimes \Lambda^{+}$, so that 
its standard  norm  is defined by 
$$|W^{+}|^{2}=\left(\frac{1}{2}\right)^{2}W^{+}_{abcd}W^{+abcd}=
W^{+}_{ABCD}W^{+ABCD},$$
where we have again used the fact that ${\varepsilon_{A}}^{A}=2$. 
For the same reason, the usual \cite{bes} 
map $W^{+}: \Lambda^{+}\to \Lambda^{+}$
is explicitly given by 
$$\psi_{ab}\mapsto \frac{1}{2}{W^{+}_{ab}}^{cd}\psi_{cd},$$
and so corresponds to  the map 
$$\psi_{AB}\mapsto {W^{+}_{AB}}^{CD}\psi_{CD}.$$
Similarly,  the usual \cite{bourg} expression $W^{+}(\psi , \psi )$ is actually 
short-hand for
$$\left(\frac{1}{2}\right)^{2}W^{+}_{abcd}\psi^{ab}\psi^{cd}=
W^{+}_{ABCD}\psi^{AB}\psi^{CD}.$$

The bundles ${\mathbb S}_{\pm}$ carry natural connections induced
by the Levi-Civita connection of $g$. These satisfy
$$\nabla_{AA'}\varepsilon_{BC}=0, ~~ 
\nabla_{AA'}\varepsilon_{B'C'}=0,$$
so the raising and lowering of spinor indices commutes with 
covariant differentiation. The curvatures of the bundles
are described by the so-called 
{\em  Ricci identities}:
 \begin{eqnarray*}
 	 \Box_{AB}\beta_{C}\hphantom{'''}& = & {W^+_{ABC}}^{D}\beta_{D}+\frac{s}{24} 
 \left(\varepsilon_{AC}\beta_{B}+\varepsilon_{BC}\beta_{A}\right) 
 \label{curv1} \\
 	\Box_{A'B'}\beta_{C}\hphantom{'} & = & -\frac{1}{2}{\ro_{A'B'C}}^{D}\beta_{D}  \\
 	\Box_{AB}\gamma_{C'} \hphantom{''}& = &  -\frac{1}{2}{\ro_{ABC'}}^{D'}\gamma_{D'}
 	\\
 	\Box_{A'B'}\gamma_{C'} & = &
 	{W^-_{A'B'C'}}^{D'}\gamma_{D'}+\frac{s}{24}  
 \left(\varepsilon_{A'C'}\gamma_{B'}+\varepsilon_{B'C'}\gamma_{A'}\right)   
 \end{eqnarray*} 
where
\begin{eqnarray*}
	\Box_{AB} \hphantom{''}& = & \nabla_{E'(A}\nabla_{B)}^{E'} , \\
	\Box_{A'B'} & = & \nabla_{E(A'}\nabla_{B')}^{E}, 
\end{eqnarray*}
so that 
$$2\nabla_{[a}\nabla_{b]}=\ve_{A'B'}\Box_{AB}+\ve_{AB}\Box_{A'B'}.$$

%
%

  \end{document}